\documentstyle[11pt,amscd,amsfonts,leqno]{amsart}

\newtheorem{Theorem}{Theorem}[section]
\newtheorem{Proposition}[Theorem]{Proposition}
\newtheorem{Lemma}[Theorem]{Lemma}
\newtheorem{Corollary}[Theorem]{Corollary}

\begin{document}

\title{Polynomial estimates, exponential curves and diophantine approximation}
\author{Dan Coman \and Evgeny A. Poletsky}
\thanks{Both authors are supported by NSF Grants.}
\date{}
\address{Department of Mathematics, Syracuse University, Syracuse,
NY 13244-1150, USA, E-mail: dcoman@@syr.edu, eapolets@@syr.edu}

\begin{abstract}
\noindent Let $\alpha\in(0,1)\setminus{\Bbb Q}$ and
$K=\{(e^z,e^{\alpha z}):\,|z|\leq1\}\subset{\Bbb C}^2$. If $P$
is a polynomial of degree $n$ in ${\Bbb C}^2$, normalized by
$\|P\|_K=1$, we obtain sharp estimates for $\|P\|_{\Delta^2}$
in terms of $n$, where $\Delta^2$ is the closed unit bidisk.
For most $\alpha$, we show that
$\sup_P\|P\|_{\Delta^2}\leq\exp(Cn^2\log n)$. However, for
$\alpha$ in a subset ${\mathcal S}$ of the Liouville numbers,
$\sup_P\|P\|_{\Delta^2}$ has bigger order of growth. We give a
precise characterization of the set ${\mathcal S}$ and study
its properties.
\end{abstract}

\maketitle

\section{Introduction}

\par The behavior of polynomials along graphs of entire transcendental
functions was recently studied in \cite{CP1,CP2,CP3}  and later
and in more general situations in \cite{Br08}. If $f$ is an
entire transcendental function and $P\in{\Bbb C}[z,w]$ is a
polynomial, the growth of the function $P(z,f(z))$ can be
estimated in terms of its uniform norm on the unit disk and the
degree of $P$. Such an estimate is called a Bernstein
inequality and it has important applications (see \cite{CP3},
\cite{Br08} and references therein). The growth estimate yields
bounds on the maximum number of zeros in a fixed disk of the
functions $P(z,f(z))$, depending only on the degree of $P$ and
$f$ \cite{CP2,CP3}. This was used in \cite{CP3} to derive
important properties of the set of algebraic numbers where the
values of $f$ are also algebraic.

\par Let $\Delta$, resp. $\Delta^2$, denote the closed unit disk,
resp. bidisk, and let ${\cal P}_n$ be the space of polynomials
$P\in{\Bbb C}[z,w]$ of degree at most $n$. The methods
introduced in \cite{CP1,CP2,CP3} involve the study of the
transcendence measures $$E_n(f)=\sup\|P\|_{\Delta^2},$$ where
$P\in{\mathcal P}_n$ is normalized by $|P(z,f(z))|\leq1$ for
$z\in\Delta$. We showed in \cite{CP3} that for any
transcendental function $f$ of finite positive order, $\log
E_n(f)$ grows like $n^2\log n$, while the maximum number of
zeros in a fixed disk of the functions $P(z,f(z))$,
$P\in{\mathcal P}_n$, grows like $n^2$, at least for an
infinite sequence of natural numbers $n$. Moreover, if $f$
verifies certain growth conditions (and in particular if $f$ is
a quasipolynomial), we proved that these estimates hold for
every $n$ (see \cite[Section 7]{CP3}).

\par It is an interesting open problem to study the behavior of
polynomials along the curve
$$\Gamma=\{(g(z),f(z)):\,z\in{\Bbb C}\}\subset{\Bbb C}^2,$$
where $g,f$ are algebraically independent entire functions. Let
$$K=\{(g(z),f(z)):\,z\in\Delta\}\subset{\Bbb C}^2.$$
Note that $K$ is pluripolar. Since the functions $g,f$ are
algebraically independent, it follows that the uniform norm
$\|\cdot\|_K$ is a norm on each vector space ${\cal P}_n$. As
${\cal P}_n$ are finite dimensional we have
$$E_n(\Gamma)=E_n(g,f):=
\sup\{\|P\|_{\Delta^2}:\;P\in{\cal P}_n,\;\|P\|_K\leq1\}<
+\infty,\;\;\forall\,n\geq0.$$

\par Once upper bounds on $E_n(\Gamma)$ are known, one can use
the classical Bernstein-Walsh inequality as in \cite{CP1} to
estimate the growth of any polynomial $P\in{\mathcal P}_n$ at
every point in terms of $n$ and $\|P\|_K$, despite the
pluripolarity of $K$:
\begin{equation}\label{e:gBW}
|P(z,w)|\leq\|P\|_KE_n(\Gamma)\exp(n\log^+\max\{|z|,|w|\}),\;
(z,w)\in{\Bbb C}^2.
\end{equation}

\par In some cases when $g,f$ have different orders of growth
certain upper bounds on $E_n(\Gamma)$ can be derived using
\cite[Theorem 2.3]{Br08}.

\par In this note we consider the simplest case of the exponential curves
$$\Gamma=\{(e^z,e^{\alpha z}):\,z\in{\Bbb C}\}\subset{\Bbb C}^2,$$
where $\alpha$ is a real irrational number. The functions $e^z$
and $e^{\alpha z}$ have the same order of growth and the same
growth of valencies. We denote in the sequel
$E_n(\alpha):=E_n(\Gamma)$. By results of Tijdeman, it is known
that, regardless of $\alpha$, the maximum number of zeros in a
fixed disk of the functions  $P(e^z,e^{\alpha z})$,
$P\in{\mathcal P}_n$, grows like $n^2$ for all $n$ (see
\cite{T}, \cite{B}).

\par We obtain here sharp estimates for $E_n(\alpha)$
and show that these estimates depend on the rate of Diophantine
approximation of $\alpha$. In contrast to the case mentioned
above when $\Gamma$ was the graph of a quasipolynomial, we see
that: 1) $E_n(\alpha)$ may grow much faster than the maximal
number of zeros in a fixed disk of the functions
$P(e^z,e^{\alpha z})$, $P\in{\mathcal P}_n$; 2) transcendental
number theory is needed to get estimates on $E_n(g,f)$ for all
$n$.
\medskip

\par We now state our results more precisely. Let
$\alpha\in(0,1)\setminus{\Bbb Q}$ and
$$e_n(\alpha)=\log E_n(\alpha).$$
Throughout the paper we denote by $p_s/q_s$ the convergents to
$\alpha$ given by its continued fractions expansion (see
Section \ref{S:Proofs}), and by $[x]$ the greatest integer
$\leq x$. We have the following:

\begin{Theorem}\label{T:enbd} Let $\alpha\in(0,1)\setminus{\Bbb Q}$
and let $p_s/q_s$, $s\geq0$, be the convergents to $\alpha$
given by its continued fractions expansion. If $q_s\leq
n<q_{s+1}$ then
$$\max\left\{\frac{n^2\log n}{2}-n^2,\,
\left[\frac{n}{q_s}\right]\log q_{s+1}-n\right\}
\leq e_n(\alpha)\leq\frac{n^2\log n}{2}+9n^2+\frac{n}{q_s}\log q_{s+1}.$$
\end{Theorem}

\par Theorem \ref{T:enbd} implies a  connection between $E_n(\alpha)$
and Diophantine approximation. Namely, $E_n(\alpha)$ provides a
lower bound for the rate of approximation of $\alpha$ by
rational numbers with denominator at most $n$.

\begin{Corollary}\label{C:dio} Let $\alpha\in(0,1)\setminus{\Bbb Q}$.
For every $n\geq1$ we have
$$\min_{1\leq k\leq n}dist(k\alpha,{\Bbb Z})\geq(2e^nE_n(\alpha))^{-1}.$$
\end{Corollary}

\par A number $\alpha\in{\Bbb R}\setminus{\Bbb Q}$ is called Diophantine
of order $\mu$, $\mu\geq2$, if there exists $\varepsilon>0$ so
that $|\alpha-p/q|>\varepsilon q^{-\mu}$, for every rational
number $p/q$. We denote by ${\mathcal D}(\mu)$ the set of such
numbers. Then
\begin{equation}\label{e:dio}
\alpha\in{\mathcal D}(\mu) \Longleftrightarrow q_{s+1}\leq Cq_s^{\mu-1},\;\forall\,s\geq0,
\end{equation}
for some constant $C>0$, where $p_s/q_s$ are the convergents to
$\alpha$ (see e.g. \cite[Appendix C]{Mil}). We let
$${\mathcal D}(\infty)=\bigcup_{\mu\geq2}{\mathcal D}(\mu)\;,\;{\mathcal L}={\Bbb R}\setminus({\Bbb Q}\cup{\mathcal D}(\infty)).$$
${\mathcal L}$ is called the set of Liouville numbers. It has
Hausdorff dimension zero (see e.g. \cite[Lemma C.4]{Mil}). By a
classical theorem of Liouville, any algebraic number of degree
$\mu$ belongs to ${\mathcal D}(\mu)$. Hence all Liouville
numbers are transcendental.

\begin{Corollary}\label{C:endio} If $\alpha\in(0,1)$ is Diophantine
of order $\mu$ then for $n\geq1$
$$\frac{n^2\log n}{2}-n^2\leq e_n(\alpha)\leq\frac{n^2\log n}{2}+9n^2+Cn,$$
where $C>0$ is a constant depending on $\alpha$.
\end{Corollary}

\par Using Theorem \ref{T:enbd}, it is in fact possible to obtain a
precise characterization of the numbers
$\alpha$ for which $e_n(\alpha)$ grows like $n^2\log n$:

\begin{Corollary}\label{C:enreg} If $\alpha\in(0,1)\setminus{\Bbb Q}$ then
$$\frac{e_n(\alpha)}{n^2\log n}=O(1)\Longleftrightarrow \frac{e_{q_s}(\alpha)}{q_s^2\log q_s}=O(1)\Longleftrightarrow \frac{\log q_{s+1}}{q_s^2\log q_s}= O(1).$$
\end{Corollary}

\par Theorem \ref{T:enbd} and its corollaries are proved in Section
\ref{S:Proofs}.  We also review there the  necessary results
about continued fractions and Diophantine approximation.

\medskip

\par Corollary \ref{C:enreg} leads us to consider the following set of irrational numbers:
$${\mathcal S}=\left\{\alpha\in(0,1)\setminus{\Bbb Q}:\,
\limsup_{s\to+\infty}\frac{\log q_{s+1}}{q_s^2\log q_s}=+\infty\right\}.$$
If $\alpha\in{\mathcal S}$ then $e_n(\alpha)$ grows faster than
$n^2\log n$ for a sequence of integers $n=q_{s_j}$, where $\log
q_{s_j+1}/(q_{s_j}^2\log q_{s_j})\to+\infty$.

\par It follows from (\ref{e:dio}) that Liouville numbers can be characterized as follows:
$$\alpha\in{\mathcal L}\Longleftrightarrow\limsup_{s\to+\infty}\frac{\log q_{s+1}}{\log q_s}=+\infty.$$
Hence ${\mathcal S}\subset{\mathcal L}$. In fact, we see from
the recursive formulas for $\{q_s\}$ (see Section
\ref{S:Proofs}) and from Theorem \ref{T:psqs} that ${\mathcal
S}$ is a ``small" subset of ${\mathcal L}$ consisting  of
transcendental numbers which are very well approximated by
rationals.

\par We study the set ${\mathcal S}$ in Section \ref{S:S}. We prove that
${\mathcal S}$ contains a dense $G_\delta$ set, hence it is
uncountable. We also prove that it has Hausdorff $h$-measure 0,
for a class of rapidly increasing functions $h$. We also
discuss the connection between ${\mathcal S}$ and certain polar
sets of Liouville numbers defined in terms of the growth of the
denominators $q_s$ given by their continued fractions
expansion.

\section{Proof of Theorem \ref{T:enbd}}\label{S:Proofs}

\par Let $\alpha\in{\Bbb R}\setminus{\Bbb Q}$. Then $\alpha$ has a unique representation as an (infinite) continued fraction
$$\alpha=[a_0;a_1,a_2,\dots]=a_0+\frac{1}{a_1+\frac{1}{a_2+\dots}}\;,$$
where all $a_j$ are integers and $a_j\geq1$ for $j\geq1$ (see
e.g. \cite[Theorem 14]{Khi}). The rational number
$$\frac{p_s}{q_s}=[a_0;a_1,a_2,\dots,a_s]=a_0+\frac{1}{a_1+\frac{1}{a_2+_{\ddots a_{s-1}+\frac{1}{a_s}}}}$$
is called the $s$-th convergent to $\alpha$. Viewing $p_s,q_s$
as polynomials in the variables $a_0,\dots,a_s$ one has  the
following recursive formulas \cite[Theorem 1]{Khi}:
$$p_s=a_sp_{s-1}+p_{s-2}\;,\;q_s=a_sq_{s-1}+q_{s-2},\;s\geq1,$$
where $p_0=a_0$, $q_0=1$, $p_{-1}=1$, $q_{-1}=0$. Moreover
\cite[Theorem 2]{Khi},
$$q_sp_{s-1}-p_sq_{s-1}=(-1)^s,$$
which implies that the fraction $p_s/q_s\in{\Bbb Q}$ is
irreducible. For $s\geq1$, $q_{s+1}>q_s$ and
$q_s\geq2^{(s-1)/2}$ \cite[Theorem 12]{Khi}. We now recall a
few properties of the convergents $p_s/q_s$, which will be
useful later.

\begin{Theorem}\label{T:psqs}\cite[Theorems 9 and 13]{Khi} For $s\geq0$,
$$(2q_{s+1})^{-1}\leq(q_{s+1}+q_s)^{-1}<\left|q_s\alpha-p_s\right|<q_{s+1}^{-1}.$$
\end{Theorem}

\par By a theorem of Lagrange, continued fractions provide the best rational approximations to $\alpha$:

\begin{Theorem}\label{T:Lagrange}\cite[Theorem 5E]{Sch80} For $s\geq0$,
$|q_s\alpha-p_s|>|q_{s+1}\alpha-p_{s+1}|$. Moreover, if
$s\geq1$, $1\leq q\leq q_s$, and if $(p,q)\neq(p_s,q_s)$,
$(p,q)\neq(p_{s-1},q_{s-1})$ then
$|q\alpha-p|>|q_{s-1}\alpha-p_{s-1}|$.
\end{Theorem}

\par Conversely, if $|d\alpha-c|>|b\alpha-a|$ for each integers $c,d$ with $1\leq d\leq b$, $c/d\neq a/b$ then $a/b$ is a convergent to $\alpha$ (\cite[Theorem 16]{Khi}).  Another result of this kind is the following theorem of Legendre:

\begin{Theorem}\label{T:Legendre}\cite[Theorem 5C]{Sch80} If $p,q$ are relatively prime integers, $q>0$ and $|q\alpha-p|<(2q)^{-1}$ then $p/q$ is a convergent to $\alpha$.
\end{Theorem}

\medskip

\par Next we develop certain estimates which will be needed in the proof of Theorem \ref{T:enbd}. Let $\alpha\in{\Bbb R}\setminus{\Bbb Q}$ and let $p_s/q_s$, $s\geq0$, be the convergents to $\alpha$ given by its continued fractions expansion. For $k\in{\Bbb N}$ we denote by $(k\alpha)$ the (unique) closest integer to $\alpha$, so
$$dist(k\alpha,{\Bbb Z})=|k\alpha-(k\alpha)|<1/2.$$

\begin{Lemma}\label{L:factorial} Let $k,x,y\in{\Bbb Z}$, $x\leq y$, $k\geq1$. Then (with $0^0:=1$)
$$\prod_{j=x}^y|j-k\alpha|\geq\left\{\begin{array}{ll}\frac{1}{2}\left(\frac{y-x}{e}\right)^{y-x},\;{\rm if}\;(k\alpha)\not\in[x,y],\\
\left(\frac{y-x}{2e}\right)^{y-x}dist(k\alpha,{\Bbb Z}),\;{\rm if}\;x\leq(k\alpha)\leq y.
\end{array}\right.$$
\end{Lemma}

\begin{pf} By Stirling's formula we have
$$e^{7/8}\leq\frac{m!}{(m/e)^m\sqrt{m}}\leq e,\;m\geq1,$$
This implies
$$\prod_{j=1}^m\left(j-\frac{1}{2}\right)=\frac{(2m)!}{2^{2m}m!}>(m/e)^m.$$

\par Let $j_0=(k\alpha)$. If $j\neq j_0$ then
$$|j-k\alpha|\geq|j-j_0|-|j_0-k\alpha|>|j-j_0|-1/2.$$
Using this we obtain for $j_0<x$,
$$\prod_{j=x}^y|j-k\alpha|\geq\prod_{j=x}^y(j-j_0-1/2)=\prod_{j=0}^{y-x}(j+x-j_0-1/2)\geq\frac{1}{2}\,(y-x)!.$$
Similarly, if $y<j_0$,
$$\prod_{j=x}^y|j-k\alpha|\geq\prod_{j=x}^y(j_0-j-1/2)=\prod_{j=0}^{y-x}(j+j_0-y-1/2)\geq\frac{1}{2}\,(y-x)!.$$

\par We assume now that $x\leq j_0\leq y$. Then, as before,
\begin{eqnarray*}
\prod_{j=x}^y|j-k\alpha|&\geq&
\prod_{j=x}^{j_0-1}(j_0-j-1/2)\,\prod_{j=j_0+1}^y(j-j_0-1/2)\;\;dist(k\alpha,{\Bbb Z})\\
&=&\prod_{j=1}^{j_0-x}(j-1/2)\,\prod_{j=1}^{y-j_0}(j-1/2)\;\;dist(k\alpha,{\Bbb Z})\\
&\geq&\left(\frac{j_0-x}{e}\right)^{j_0-x}\left(\frac{y-j_0}{e}\right)^{y-j_0} dist(k\alpha,{\Bbb Z}).
\end{eqnarray*}
The function $f(t)=(t-x)\log(t-x)+(y-t)\log(y-t)$ attains its
minimum on the interval $[x,y]$ at $t=(x+y)/2$, so
$$f(t)\geq(y-x)\log\left(\frac{y-x}{2}\right).$$
This implies
$$\prod_{j=x}^y|j-k\alpha|\geq\left(\frac{y-x}{2e}\right)^{y-x}dist(k\alpha,{\Bbb Z}).$$
\end{pf}

\par The following result provides lower estimates for the function
$$D_\alpha(n)=\prod_{k=1}^n dist(k\alpha,{\Bbb Z}).$$

\begin{Lemma}\label{L:dalpha} If $q_s\leq n<q_{s+1}$ then
$D_\alpha(n)\geq(2n)^{-n}q_{s+1}^{-n/q_s}$.
\end{Lemma}

\begin{pf} We consider the sets
$$S_j=\left\{k\in{\Bbb N}:\,k\leq n,\;\frac{(k\alpha)}{k}=\frac{p_j}{q_j}\right\},\;0\leq j\leq s,\;\;
S_{s+1}=([1,n]\cap{\Bbb N})\setminus\bigcup_{j=0}^sS_j.$$
For $1\leq k\leq n$, suppose that
$$dist(k\alpha,{\Bbb Z})=|k\alpha-(k\alpha)|<(2k)^{-1}.$$
Theorem \ref{T:Legendre} implies that $(k\alpha)/k=p_j/q_j$ for
some $j\leq s$, so $k\in S_j$. We conclude that for $k\in
S_{s+1}$
$$dist(k\alpha,{\Bbb Z})\geq(2k)^{-1}\geq(2n)^{-1}.$$
Hence
$$\prod_{k\in S_{s+1}}dist(k\alpha,{\Bbb Z})\geq(2n)^{-|S_{s+1}|}.$$

\par Since $p_j/q_j$ is irreducible it follows that the sets $S_j$, $j\leq s$, are disjoint and
$$dist(k\alpha,{\Bbb Z})=|k\alpha-(k\alpha)|\geq|q_j\alpha-p_j|\geq(2q_{j+1})^{-1}, \;k\in S_j.$$
Here the last inequality follows by Theorem \ref{T:psqs}.
Moreover, if $k\in S_s$ then $q_s|k$, so $|S_s|\leq n/q_s$.
Hence
\begin{eqnarray*}
&&\prod_{k\in S_j}dist(k\alpha,{\Bbb Z})\geq(2q_{j+1})^{-|S_j|}\geq(2n)^{-|S_j|},\;0\leq j<s,\\
&&\prod_{k\in S_s}dist(k\alpha,{\Bbb Z})\geq(2q_{s+1})^{-|S_s|}\geq2^{-|S_s|}q_{s+1}^{-n/q_s}.
\end{eqnarray*}

\par Note that $|S_0|+\dots+|S_{s+1}|=n$. We conclude that
$$D_\alpha(n)=\prod_{j=0}^{s+1}\prod_{k\in S_j}dist(k\alpha,{\Bbb Z})\geq(2n)^{-n}q_{s+1}^{-n/q_s}.$$
\end{pf}

\begin{Lemma}\label{L:dalpha2} If $q_s\leq n<q_{s+1}$ and $0\leq m\leq n$ then
$$D_\alpha(m)D_\alpha(n-m)\geq2^{-n}n^{-2n}q_{s+1}^{-n/q_s}.$$
\end{Lemma}

\begin{pf} There exist integers $j,l$ so that $q_j\leq m<q_{j+1}$ and $q_l\leq n-m<q_{l+1}$. Note that $m^m(n-m)^{(n-m)}\leq n^n$, so by Lemma \ref{L:dalpha},
$$D_\alpha(m)D_\alpha(n-m)\geq(2n)^{-n}q_{j+1}^{-m/q_j}q_{l+1}^{-(n-m)/q_l}.$$
If $\max\{j,l\}<s$ then
$$q_{j+1}^{-m/q_j}q_{l+1}^{-(n-m)/q_l}\geq q_s^{-m/q_j-(n-m)/q_l}\geq n^{-n}.$$
If $l=s>j$ then
$$q_{j+1}^{-m/q_j}q_{l+1}^{-(n-m)/q_l}\geq n^{-n}q_{s+1}^{-n/q_s}.$$
Finally, if $j=l=s$ then
$$q_{j+1}^{-m/q_j}q_{l+1}^{-(n-m)/q_l}=q_{s+1}^{-n/q_s}.$$
\end{pf}

\vspace{5mm}

\par\noindent {\em Proof of Theorem \ref{T:enbd}.} Recall that ${\dim}\,{\cal
P}_n=N+1$, where $N=(n^2+3n)/2$.

\par We start by proving the upper bound for $e_n(\alpha)$. Let us introduce the following notation. For any polynomial
$R(\lambda)=\sum_{j=0}^mc_j\lambda^j$ we denote by $D_R$ the
constant-coefficient differential operator
$$D_R=R\left(\frac{d}{dz}\right)=\sum_{j=0}^mc_j\frac{d^j}{dz^j}\,.$$
Then for any integer $t\geq0$ and any $a\in{\Bbb C}$ we have
\begin{equation}\label{e:diff}
\left.D_R[z^te^{a z}]\,\right|_{z=0}=\sum_{j\geq t}c_j\frac{j!}{(j-t)!}\,a^{j-t}=\left.\frac{d^tR}{d\lambda^t}\,
\right|_{\lambda=a}=R^{(t)}(a).
\end{equation}
Fix now $P\in{\cal P}_n$, $n\geq1$, with $\|P\|_K\leq1$. We
write
$$
P(z,w)=\sum_{j+k\leq n}c_{jk}z^jw^k,\;f(z):=P(e^z,e^{\alpha z})=\sum_{j+k\leq n}c_{jk}e^{(j+k\alpha)z}.
$$

\par We will estimate the coefficients $c_{lm}$ of $P$ by using the differential operators given by the
polynomials of degree $N$,
$$R_{lm}(\lambda)=\prod_{j+k\leq n,(j,k)\neq(l,m)}(\lambda-j-k\alpha)=\sum_{t=0}^Na_t\lambda^t.$$
Since the coefficients $a_t$ are elementary symmetric functions
of the roots of $R_{lm}$ it follows that for $\lambda\geq0$
$$\sum_{t=0}^N|a_t|\lambda^t\leq\prod_{j+k\leq n,(j,k)\neq(l,m)}(\lambda+|j+k\alpha|)\leq(\lambda+n)^N,$$
where for the last inequality we used $|j+k\alpha|\leq j+k\leq
n$, since $0<\alpha<1$.

\par By (\ref{e:diff}) we have
$$
D_{R_{lm}}f(z)\mid_{z=0}=c_{lm}\beta_{lm}\;,\;\;\beta_{lm}=\prod_{j+k\leq n,(j,k)\neq(l,m)}(l-j+(m-k)\alpha).
$$
By Cauchy's estimates $|f^{(t)}(0)|\leq t!\leq N^t$ for $t\leq
N$, so we obtain
$$|\left.D_{R_{lm}}f(z)\right|_{z=0}|=\left|\sum_{t=0}^Na_tf^{(t)}(0)\right|\leq\sum_{t=0}^N|a_t|N^t\leq
(N+n)^N.$$
Therefore
\begin{equation}\label{e:cbeta}
\log(|c_{lm}\beta_{lm}|)\leq N\log(N+n)\leq n^2\log n+3.7n^2.
\end{equation}

\par Next we obtain lower estimates on $|\beta_{lm}|$. We have
$$|\beta_{lm}|\geq\prod_{k=0,k\neq m}^n\prod_{j=0}^{n-k}|l-j+(m-k)\alpha|=A_1A_2,$$
where
$$A_1=\prod_{k=0}^{m-1}\prod_{j=0}^{n-k}|j-l-(m-k)\alpha|=\prod_{k=1}^m\prod_{j=-l}^{n-m-l+k}|j-k\alpha|,$$
$$A_2=\prod_{k=m+1}^n\prod_{j=0}^{n-k}|l-j-(k-m)\alpha|=\prod_{k=1}^{n-m}\prod_{j=l+m-n+k}^l|j-k\alpha|.$$
By Lemma \ref{L:factorial}
\begin{eqnarray*}
A_1&\geq&D_\alpha(m)\prod_{k=1}^m\left(\frac{n-m+k}{2e}\right)^{n-m+k},\\
A_2&\geq&D_\alpha(n-m)\prod_{k=1}^{n-m}\left(\frac{n-m-k}{2e}\right)^{n-m-k}.
\end{eqnarray*}
Thus, using Lemma \ref{L:dalpha2},
\begin{eqnarray*}
|\beta_{lm}|&\geq&D_\alpha(m)D_\alpha(n-m)\prod_{k=n-m+1}^n\left(\frac{k}{2e}\right)^k\times
\prod_{k=0}^{n-m-1}\left(\frac{k}{2e}\right)^k\\
&\geq&2^{-n}n^{-2n}q_{s+1}^{-n/q_s}
\left(\frac{2e}{n-m}\right)^{n-m}\prod_{k=1}^n\left(\frac{k}{2e}\right)^k\\
&\geq&2^{-n}n^{-3n}(2e)^{-n^2}q_{s+1}^{-n/q_s}\prod_{k=1}^nk^k.
\end{eqnarray*}
We have (see e.g. \cite[Lemma 2.1]{CP1})
$$\sum_{k=1}^n k\log k\geq\frac{n^2\log n}{2}-\frac{n^2}{4}\;.$$
This yields
$$\log|\beta_{lm}|\geq\frac{n^2\log n}{2}-4.2n^2-\frac{n}{q_s}\log q_{s+1}.$$
Using (\ref{e:cbeta}) we obtain
$$\log|c_{lm}|\leq\frac{n^2\log n}{2}+7.9n^2+\frac{n}{q_s}\log q_{s+1}.$$
Since $\|P\|_{\Delta^2}\leq\sum|c_{jk}|\leq(N+1)\max|c_{jk}|$,
we conclude that
$$e_n(\alpha)\leq\frac{n^2\log n}{2}+9n^2+\frac{n}{q_s}\log q_{s+1}.$$

\medskip

We now proceed to prove the lower bound for $e_n(\alpha)$.
There exists a non-trivial polynomial $P\in{\mathcal P}_n$ so
that the function $P(e^z,e^{\alpha z})$ has a zero of order at
least $N=(n^2+3n)/2$ at 0. Using (\ref{e:gBW}) and repeating
the argument in the proof of \cite[Proposition 1.3]{CP1}, we
obtain that
$$e_n(\alpha)\geq N\log n-n^2.$$

\par Consider the polynomial $P(z,w)=(z^{p_s}-w^{q_s})^{[n/q_s]}$. Since $0<\alpha<1$, we have $0\leq p_s\leq q_s$ for every $s$, so $P\in{\mathcal P}_n$. Note that $\|P\|_{\Delta^2}=2^{[n/q_s]}$. If $|z|\leq1$ we have by Theorem \ref{T:psqs}
$$|(q_s\alpha-p_s)z|\leq q_{s+1}^{-1}\leq1.$$
Using that $|1-e^\zeta|\leq2|\zeta|$ for $|\zeta|\leq1$, we
obtain
$$|P(e^z,e^{\alpha z})|\leq\left|e^{p_sz}\left(1-e^{(q_s\alpha-p_s)z}\right)\right|^{[n/q_s]}\leq e^n(2q_{s+1}^{-1})^{[n/q_s]},\;|z|\leq1.$$
Therefore
$$E_n(\alpha)\geq\|P\|_{\Delta^2}/\|P\|_K\geq q_{s+1}^{[n/q_s]}e^{-n},$$
and the proof is complete. $\Box$

\vspace{5mm}

\par\noindent {\em Proof of Corollary \ref{C:dio}.} Theorems \ref{T:Lagrange} and \ref{T:psqs} show that if $q_s\leq n<q_{s+1}$ then
$$\min_{1\leq k\leq n}dist(k\alpha,{\Bbb Z})=|q_s\alpha-p_s|\geq1/(2q_{s+1}),$$
while the lower bound for $e_n(\alpha)$ from Theorem
\ref{T:enbd} implies $\log q_{s+1}\leq e_n(\alpha)+n$. It
follows that for all $n\geq 1$ we have
$$\min_{1\leq k\leq n}dist(k\alpha,{\Bbb Z})\geq(2e^nE_n(\alpha))^{-1}.\;\;\;\Box$$

\vspace{5mm}

\par\noindent {\em Proof of Corollary \ref{C:endio}.} The upper estimate follows immediately from Theorem \ref{T:enbd}, since by (\ref{e:dio})
$$\frac{\log q_{s+1}}{q_s}\leq\frac{\log C}{q_s}+(\mu-1)\frac{\log q_s}{q_s}
\leq\log C+\frac{\mu-1}{2}\;.\;\;\;\Box$$

\vspace{5mm}

\par\noindent {\em Proof of Corollary \ref{C:enreg}.} Assume first that $e_{q_s}(\alpha)\leq Cq_s^2\log q_s$ for all $s$, where $C$ is a constant. By the lower estimate in Theorem \ref{T:enbd} applied for $n=q_s$, we get
$$\log q_{s+1}\leq e_{q_s}(\alpha)+q_s\leq(C+1)q_s^2\log q_s.$$

\par Assume now that $\log q_{s+1}\leq Cq_s^2\log q_s$ for all $s$, where $C$ is a constant. Given $n$, there is a unique $s$ so that $q_s\leq n<q_{s+1}$. By Theorem \ref{T:enbd},
$$e_n(\alpha)\leq\frac{n^2\log n}{2}+9n^2+\frac{n}{q_s}\log q_{s+1}\leq(C+10)n^2\log n.\;\;\;\Box$$

\section{The set ${\mathcal S}$}\label{S:S}

\par Let $E\subset{\Bbb C}$ and $h(r)$, $0\leq r\leq r_0$, be a continuous increasing function with $h(0)=0$. Given
$\delta>0$ we define
$${\mathcal H}^h_\delta(E)=\inf\sum_{n=1}^\infty h({\rm diam}\,A_n/2),$$
where the infimum is taken over all coverings $\{A_n\}$ of $E$
with bounded sets $A_n$ of diameter less than $\delta$. As
$\delta\searrow0$ the quantities ${\mathcal H}^h_\delta(E)$
increase, so the limit
$${\mathcal H}^h(E)=\lim_{\delta\to0}{\mathcal H}^h_\delta(E)$$
exists and is called the Hausdorff $h$-measure of $E$ (see e.g.
\cite[p. 196]{L}). We recall that if $h(r)=1/\log(1/r)$ then
${\mathcal H}^h$ is called the {\em logarithmic measure}. A set
$E\subset{\Bbb C}$ of finite logarithmic measure is polar
\cite[Theorem 3.14]{L}.

\par We assume now that $h(r)$, $0\leq r\leq r_0$, is a continuous increasing function so that
$$\sum_{n=N}^\infty n\,h\left(n^{-n^2}\right)<+\infty.$$
An example of such a function is
$$h(r)=\frac{1}{\log\frac{1}{r}\left(\log\log\log\frac{1}{r}\right)^p}\;,\;\;p>1.$$

\begin{Proposition}\label{P:gdh} If $h$ is as above then ${\mathcal H}^h(S)=0$. Moreover, ${\mathcal S}$ contains a dense $G_\delta$ set, hence it is uncountable.
\end{Proposition}

\begin{pf} Note that by Theorem \ref{T:psqs} and the definition of ${\mathcal S}$ we have the following: if
$\alpha\in{\mathcal S}$ then there exist infinitely many
rational numbers $p_s/q_s$ so that
$$|\alpha-p_s/q_s|<q_{s+1}^{-1}<q_s^{-q_s^2}.$$
Let $r(n)=n^{-n^2}$ and define
$$A_n=\bigcup_{m=1}^n\left(\frac{m}{n}-r(n),\frac{m}{n}+r(n)\right).$$
It follows that
$${\mathcal S}\subset\limsup A_n=\bigcap_{k=1}^\infty\bigcup_{n\geq k}A_n.$$

\par Fix $\delta>0$. If $k$ is large enough so that $2r(k)<\delta$, then by the definition of
${\mathcal H}^h_\delta$
$${\mathcal H}^h_\delta({\mathcal S})\leq{\mathcal H}^h_\delta\left(\cup_{n\geq k}A_n\right)\leq
\sum_{n\geq k}n\,h(r(n)).$$
Since $\sum_{n\geq1}n\,h(r(n))<+\infty$, it follows that
${\mathcal H}^h_\delta({\mathcal S})=0$ for all $\delta>0$, so
${\mathcal H}^h(S)=0$.

\par We now let $r'(n)=e^{-n^3}$ and define
$$A'_n=\bigcup_{m=1,(m,n)=1}^n\left(\frac{m}{n}-r'(n),\frac{m}{n}+r'(n)\right),\;
G=\limsup A'_n=\bigcap_{k=1}^\infty\bigcup_{n\geq k}A'_n.$$
Here $(m,n)$ denotes the greatest common divisor of $m,n$. By
Baire's theorem, $G$ is a dense $G_\delta$ set and hence it is
uncountable.

\par Let us show that $G\subset{\mathcal S}$. If $\alpha\in G$ there exists a sequence of rational numbers $m_k/n_k$ with $(m_k,n_k)=1$ and $n_k\to+\infty$, so that $|\alpha-m_k/n_k|<r'(n_k)$. Thus
$$|n_k\alpha-m_k|<n_ke^{-n_k^3}<(2n_k)^{-1}.$$
This implies that $\alpha$ is irrational. Indeed, if
$\alpha=p/q\in{\Bbb Q}$ with $(p,q)=1$ then for $n_k>q$ we have
$$q^{-1}\leq|n_k\alpha-m_k|<n_ke^{-n_k^3},$$
which yields a contradiction.

\par Since $|n_k\alpha-m_k|<(2n_k)^{-1}$ we see by Theorem \ref{T:Legendre} that $m_k/n_k$ is a convergent to $\alpha$, so $m_k=p_s$ and  $n_k=q_s$ for some $s$. Using Theorem \ref{T:psqs} we obtain
$$(2q_{s+1})^{-1}<|q_s\alpha-p_s|<q_se^{-q_s^3}\Longrightarrow\frac{\log q_{s+1}}{q_s^2\log q_s}>\frac{q_s}{\log q_s}-o(1).$$
We conclude that $\alpha\in{\mathcal S}$.
\end{pf}

\medskip

\par\noindent{\bf Remark.} An argument similar to the one used to prove ${\mathcal H}^h({\mathcal S})=0$ shows that the above dense $G_\delta$ set $G$ has zero logarithmic measure, hence it is polar.

\medskip

\par We conclude this section by considering certain polar sets of irrational numbers related to ${\mathcal S}$. Given a sequence
$\varepsilon:{\Bbb N}\to(0,+\infty)$ we introduce the sets
\begin{eqnarray*}
{\mathcal T}(\varepsilon)&=&\{\alpha\in(0,1)\setminus{\Bbb Q}:\,\limsup_{s\to\infty}\varepsilon(q_s)\log q_{s+1}=+\infty\},\\
{\mathcal U}(\varepsilon)&=&\{\alpha\in(0,1)\setminus{\Bbb Q}:\,\limsup_{n\to\infty}\varepsilon(n)e_n(\alpha)=+\infty\}.
\end{eqnarray*}
Our interest will be in sequences $\varepsilon$ that in some
sense decrease rapidly to 0. We have the following:

\begin{Proposition}\label{P:polar} (i) If $\varepsilon$ satisfies
$\displaystyle\sum_{n=1}^\infty n\,\varepsilon(n)<\+\infty$
then the set ${\mathcal T}(\varepsilon)$ is polar.

\par (ii) If $\varepsilon$ is given by
$\varepsilon(n)=\left(x(n)n^2\log n\right)^{-1}$, $n\geq1$,
where $x(n)\geq1$ is an increasing sequence, then   ${\mathcal
T}(\varepsilon)={\mathcal U}(\varepsilon)\subset{\mathcal S}$.
\end{Proposition}

\begin{pf} $(i)$ Consider the function
$$v(\zeta)=\sum_{n=1}^\infty\varepsilon(n)\sum_{m=1}^n\log\frac{|\zeta-m/n|}{3}\;,\;\;|\zeta|<2.$$
We have that
$$v(i)\geq-\log3\,\sum_{n=1}^\infty n\,\varepsilon(n)>-\infty,$$
so $v$ is a negative subharmonic function in $\{|\zeta|<2\}$.
If $\alpha\in(0,1)\setminus{\Bbb Q}$ it follows from Theorem
\ref{T:psqs} that  $|\alpha-p_s/q_s|<q_{s+1}^{-1}$, so
$v(\alpha)<-\varepsilon(q_s)\log q_{s+1}$, for every $s$. Hence
if $\alpha\in{\mathcal T}(\varepsilon)$ then
$v(\alpha)=-\infty$.

\medskip

\par $(ii)$ Clearly ${\mathcal T}(\varepsilon)\subset{\mathcal S}$. Using the lower bound for $e_n(\alpha)$ from Theorem \ref{T:enbd} with $n=q_s$ we obtain
$$\varepsilon(q_s)e_{q_s}(\alpha)\geq\varepsilon(q_s)\log q_{s+1}-\varepsilon(q_s)q_s.$$
Since $\varepsilon(q_s)q_s\to0$, we deduce that ${\mathcal
T}(\varepsilon)\subset{\mathcal U}(\varepsilon)$.

\par Conversely, if $\alpha\in{\mathcal U}(\varepsilon)$ there exists a sequence of integers $n_j\to+\infty$ so that $\varepsilon(n_j)e_{n_j}(\alpha)\to+\infty$. We have $q_{s_j}\leq n_j<q_{s_j+1}$ for a unique $s_j$, so by Theorem \ref{T:enbd},
\begin{eqnarray*}
\varepsilon(n_j)e_{n_j}(\alpha)&\leq&\frac{10}{x(n_j)}+\frac{\log q_{s_j+1}}{q_{s_j}x(n_j)n_j\log n_j}\\
&\leq&
10+\frac{\log q_{s_j+1}}{x(q_{s_j})q_{s_j}^2\log q_{s_j}}=10+\varepsilon(q_{s_j})\log q_{s_j+1}.
\end{eqnarray*}
We conclude that $\alpha\in{\mathcal T}(\varepsilon)$.
\end{pf}

\medskip

\par\noindent{\bf Remark.} There exists a sequence $\varepsilon$ which verifies the hypothesis of Proposition \ref{P:polar} $(i)$ for which ${\mathcal U}(\varepsilon)=(0,1)\setminus{\Bbb Q}$. Indeed, we let
$$\varepsilon(n)=\left\{\begin{array}{ll}
n^{-3},\;\;\;{\rm if}\;n\in{\Bbb N}\setminus\left\{2^k:\,k\in{\Bbb N}\right\},\\
n^{-2},\;\;\;{\rm if}\;n\in\left\{2^k:\,k\in{\Bbb N}\right\}.
\end{array}\right.$$
Clearly $\sum_{n=1}^\infty n\,\varepsilon(n)<\+\infty$. Let
$\alpha\in(0,1)\setminus{\Bbb Q}$. By Theorem \ref{T:enbd} we
obtain
$$\varepsilon(n)e_{n}(\alpha)\geq\varepsilon(n)\frac{n^2\log n}{2}-\varepsilon(n)n^2=\frac{\log n}{2}-1,\;{\rm if}\;n=2^k,\;k\in{\Bbb N}.$$

\end{document}